\input amssym.def
\input amssym
\magnification=1200
\parindent0pt
\hsize=16 true cm
\baselineskip=13  pt plus .2pt

\def\Z{\Bbb Z}
\def\A{\Bbb A}
\def\S{\Bbb S}
\def\C{\Bbb C}

\def\R{\Bbb R}

\centerline {\bf On finite simple and nonsolvable groups acting on closed
4-manifolds}

\bigskip

\centerline {Mattia Mecchia and Bruno Zimmermann}

\medskip

\centerline {Universit\`a degli Studi di Trieste}
\centerline {Dipartimento di Matematica e Informatica}
\centerline {34100 Trieste, Italy}
\centerline {mecchia@dmi.units.it,  zimmer@units.it}

\bigskip \bigskip

Abstract.  {\sl  We show that the only finite nonabelian simple groups which
admit a locally linear, homologically trivial action on a closed simply
connected 4-manifold $M$ (or on a 4-manifold with trivial first homology)  are
the alternating groups $\A_5$, $\A_6$ and the linear fractional group ${\rm
PSL(2,7)}$ (we note that for homologically nontrivial actions all finite groups
occur). The situation depends strongly on the second Betti number $b_2(M)$ of
$M$ and has been known before if $b_2(M)$ is different from two, so the main
new result of the paper concerns the case $b_2(M)=2$. We prove that the only
simple group that occurs in this case is $\A_5$, and then give  a short list of
finite nonsolvable groups which contains all candidates for actions of such
groups.}

\bigskip
\bigskip

{\bf 1. Introduction}

\medskip

We are interested in actions of finite groups on closed orientable 4-manifolds.
All actions in the present paper will be {\it locally linear}, {\it faithful}
and {\it orientation-preserving}. An action is  locally linear if the isotropy
group of  each  point leaves invariant a neighbourhood of the point which is
equivariantly homeomorphic to an invariant neighbourhood of the origin in a
linear action on some Euclidean space $\R^n$ (e.g. smooth actions).

\medskip

It has been shown in [MZ1] that the only finite nonabelian simple groups acting
on a homology 4-sphere are the alternating groups $\A_5$ and $\A_6$, and from
this a short list of finite nonsolvable groups is deduced containing all
candidates for an action on a homology 4-sphere (and in particular on the
4-sphere; the corresponding situation in dimension three is considered in
[MZ2,3] and [Z2]).  On the other hand, since each finitely presented group is
the fundamental group of a closed 4-manifold, each finite group $G$ admits a
free action on a simply connected closed 4-manifold (the universal covering of
a closed 4-manifold with fundamental group $G$); as a consequence of the
Lefschetz fixed point Theorem, such a free action has to act nontrivially on
homology. For homologically trivial actions on simply connected 4-manifolds,
and more generally on 4-manifolds $M$ with trivial first homology $H_1(M)$,
there are again strong restrictions. Building on previous work of various
authors and concentrating on the basic case of nonabelian simple groups, the
following holds.

\bigskip

{\bf Theorem.}  {\sl  Let $G$ be a finite  nonabelian simple group which admits
a homologically trivial action on a closed 4-manifold $M$ with trivial first
homology. Then $G$ is isomorphic to $\A_5$, $\A_6$ or ${\rm PSL(2,7)}$.}

\bigskip

For various cases the Theorem has been known previously (and also in greater
generality);  in fact,  the situation depends strongly on the second Betti
number $b_2(M)$ of $M$, and we shall discuss the different cases in the
following; $M$ will always denote a closed 4-manifold with trivial first
homology $H_1(M)$.

\bigskip

{\bf  I.  The case $b_2(M) \ge 3$}

\medskip

In this case the possible finite groups which admit an action are very restricted;
in particular, no nonabelian simple groups occur. The following is the main result
of [Mc1].
\bigskip

{\bf Theorem 1.} ([Mc1])  {\sl Let $G$ be a finite group with a homologically
trivial action on a closed 4-manifold $M$ with trivial first homology.

i) If $b_2(M) \ge 3$ then $G$ is  abelian of rank at most two (cyclic or a
product of two cyclic groups), and $G$ has a global fixed point.

ii) If $b_2(M) \ge 2$ and $G$ has a global fixed point then again $G$ is
abelian of rank at most two.}

\bigskip

Also in the next case a complete classification is still known.

\bigskip

{\bf II.  The case $b_2(M)=1$}

\medskip

Now $M$ is a homology complex projective plane $\C P^2$, so a reference model
here is the group of projectivities  PGL$(3,\C)$ of $\C P^2$. The case has been
considered in [W1] and [HL] for homologically trivial actions of arbitrary
finite groups (see also [W2]), and the following holds.
\bigskip

{\bf Theorem 2.} ([W1],[HL]) {\sl  Let $G$ be a finite group which admits a
homologically trivial action on a a closed 4-manifold $M$ with $b_2(M)=1$ and
trivial first homology (e.g. the complex projective plane $\C P^2$). Then $G$
is isomorphic to a subgroup of ${\rm PGL}(3,\C)$ and in particular, if $G$ is
nonabelian simple, to $\A_5$, $\A_6$ or ${\rm PSL(2,7)}$.}

\bigskip

These are exactly the finite nonabelian simple subgroups of  PGL$(3,\C)$. A
main ingredient of the proof of Theorem 2  is the classification of the finite
simple groups of 2-rank at most two (i.e.  without subgroups $(\Z_2)^3$).
\bigskip

In the two remaining cases, a complete classification seems still far, so we
concentrate on the basic case of finite simple groups in the following.

\bigskip
\vfill \eject

{\bf III.  The case $b_2(M)=0$ of a homology 4-sphere}

\medskip

Now the following holds.

\bigskip

{\bf Theorem 3.}  ([MZ1]) {\sl  A finite nonabelian simple group acting on a
homology 4-sphere, and in particular on the 4-sphere $S^4$, is isomorphic to an
alternating group $\A_5$ or $\A_6$.}

\bigskip

This result is used in [MZ1] to obtain a short list of finite nonsolvable groups
which contains  all candidates for actions on a homology 4-sphere of such groups.  A
reference model here is the orthogonal group SO(5) acting on $S^4$.  We note that
now subgroups $(\Z_2)^4$ may occur, so the 2-rank of a finite simple group  acting
on a homology 4-sphere might, in principle, be equal to  four; in fact, the main
ingredient of the proof is the Gorenstein-Harada classification of the  finite
simple groups of sectional 2-rank at most four (i.e. each 2-subgroup is generated by
at most four elements).

\bigskip

Finally, we come to the last and main case of the present paper.

\bigskip

{\bf IV.  The case $b_2(M) = 2$}

\medskip

The main results of the present paper are the following  two theorems.

\bigskip

{\bf Theorem 4.}  {\sl Let $G$ be a finite nonabelian simple group which admits
an action on a closed 4-manifold $M$ with $b_2(M)=2$  and trivial first
homology (e.g. $S^2 \times S^2$). Then $G$ is isomorphic to the alternating
group $\A_5$, and $M$ has the intersection form $\pmatrix { 0  &  1 \cr  1  & 0
\cr}$ of $S^2 \times S^2$.}

\bigskip

In particular, if $M$ is simply connected then by [F] it is homeomorphic to
$S^2 \times S^2$, so a reference model here is $S^2 \times S^2$ and the  group
${\rm SO}(3) \times {\rm SO}(3)$  that consists of isometries  of $S^2 \times
S^2$ that act trivially on homology  (whose only finite simple subgroup is
$\A_5$). Again subgroups $(\Z_2)^4$ might in principle occur and, as for
Theorem 3, a main tool of the proof will be the Gorenstein-Harada
classification of the finite simple groups of sectional 2-rank at most four (as
well as part ii) of Theorem 1).

\smallskip

On the basis of Theorem 4 we obtain a short list of finite nonsolvable groups
which contains  all candidates for actions on this class of manifolds; the
result is summarized in the following theorem.

\bigskip

{\bf  Theorem 5.}  {\sl Let $G$ be a finite nonsolvable group that   admits  a
homologically trivial action on a closed 4-manifold M with $b_2(M)=2$ and
trivial first homology. Then  $G$ contains, of index at most two,  a normal
subgroup isomorphic to one of the following groups:
$$\Bbb{A}_5\times C,\;\;
\Bbb{A}_5^*\times_{\Bbb{Z}_2}C,\;\;  \Bbb{A}_5\times \Bbb{A}_5,\;\;  \Bbb{A}_5\times
\Bbb{A}_4$$   where $C$ is a cyclic group.}

\bigskip

This is close to the list of the finite nonsolvable subgroups of ${\rm SO}(3) \times
{\rm SO}(3)$, except that we are not able to exclude the binary dihedral group
$\Bbb{A}_5^*$ at moment (we suppose that it does not act).  Some information about
the possible 2-extensions can be deduced from the proof of Theorem 5.

\bigskip

{\bf 2. Proof of Theorem 4}

\medskip

In the following, $M$ will always denote a closed 4-manifold with  $b_2(M)=2$
and trivial first homology $H_1(M)$,  and $G$ a finite nonabelian simple group
acting faithfully and locally linearly on $M$.  Since the finite subgroups of
GL$(2,\Z)$ are cyclic or dihedral and $G$ is nonabelian simple, the action of
$G$ is homologically trivial. We start with some preliminary results.

\bigskip

{\bf Lemma 1.} {\sl  Let $g$ be an orientation-preserving periodic map of $M$ which
is not the identity and acts trivially on the homology of $M$. Then the fixed point
set of $g$ is of one of the following types:

i) four isolated points;

ii) a 2-sphere $S^2$ and two isolated points;

iii) two 2-spheres $S^2$.}

\medskip

{\it Proof.}  By a version of the Lefschetz fixed point Theorem (see e.g. [TD]), the
Euler characteristic of the fixed point set of $g$ equals the alternating sum of the
traces of the maps induced by $g$ on the rational homology $H_*(M;\Bbb Q)$ of $M$;
since $g$ acts trivially on homology and $H_1(M)=0$ this is equal to four. By [E,
Prop.2.4], the fixed point set of $g$ has no 1-dimensional components and  consists
of isolated points and 2-spheres; this leaves the three possibilities of the Lemma,
thus finishing its proof.
\bigskip

{\bf Lemma 2.} {\sl Let $S$ be a finite 2-group which admits a faithful,
homologically trivial action on $M$. Then $S$ is generated by at most four
elements. In particular, $G$ has sectional 2-rank at most four (each 2-subgroup
is generated by at most four elements).}
\medskip

{\it Proof.}  Let $g$ be a central involution in $S$. The possible fixed point sets
of $g$ are listed in Lemma 1.  If the fixed point set Fix($g$) consists of a
2-sphere $S^2$ and two isolated points then the 2-sphere $S^2$ is invariant under
$S$. By  a result of Edmonds [E], see also [Mc1, Theorem 2], the action on $S^2$
is  orientation preserving. The subgroup of $S$ acting trivially on $S^2$ is cyclic
(since $S$ acts homologically trivial and hence orientation-preserving on $M$), its
factor group acts faithfully on $S^2$ and, being a 2-group, is a subgroup of a
dihedral group. Clearly  $S$ is generated by at most three elements.
\medskip

Suppose that Fix($g$) consists of two 2-spheres. Then a subgroup of index at most
two of $S$ leaves invariant both 2-spheres and, considering the first case, $S$ is
generated by at most four elements.
\medskip

Finally, suppose  that Fix($g$) consists of four isolated points, invariant under
$S$. Let $S_0$ be the subgroup of $S$ fixing one of these four points, of index at
most four.  By  Theorem 1 ii), $S_0$ is abelian of rank at most two and hence $S$ is
generated by at most four elements.
\medskip

This finishes the proof of Lemma 2.

\bigskip

We apply the Gorenstein-Harada classification of the finite simple groups of
sectional 2-rank at most four (see [G1,p.6], [Su2,chapter 6, Theorem 8.12]). By
Lemma 2,  $G$ has sectional 2-rank at most four and hence is one of the groups in
the Gorenstein-Harada list; the groups are the following ($q$ denotes an odd prime
power):

$${\rm PSL}(m,q), \;\; {\rm PSU}(m,q), \; \;  \; m \le 5,$$ $${\rm G}_2(q),\;\;
^3{\rm D}_4(q),\;\;  {\rm PSp}(4,q), \; \; ^2{\rm G}_2(3^{2m+1}) \;\; (m \ge 1),$$
$${\rm PSL}(2,8),\;\; {\rm PSL(2,16)},\;\; {\rm PSL}(3,4),\;\; {\rm PSU}(3,4),\;\;
{\rm Sz}(8),$$
$$\A_m \;\; (7 \le m \le 11), \;\; {\rm M}_i \;\;(i \le 23),\;\; {\rm J}_i  \;\;
(i\le 3),  \;\; {\rm McL}, \;\; {\rm Ly}.$$

\bigskip

In the following, we will exclude all of these groups except $\A_5$. We
consider first the linear fractional groups ${\rm PSL(2,p)}$, for a prime $p
\ge 5$. The group ${\rm PSL(2,p)}$ has a metacyclic subgroup (semidirect
product) $\Z_p \ltimes \Z_{(p-1)/2}$ (represented by all upper triangular
matrices), with an effective action of $\Z_{(p-1)/2}$ (the diagonal matrices)
on the normal subgroup $\Z_p$ (the matrices having both entries one on the
diagonal).

\bigskip

{\bf Lemma 3.} {\sl For an odd prime $p$ and an integer $q \ge 2$, let $U = \Z_p
\ltimes \Z_q$ be a metacyclic group, with an effective action of $\Z_q$ on the
normal subgroup $\Z_p$.

i) If $U$ admits a faithful, orientation-preserving action
on a homology 3-sphere then  $q=2$.

ii) If  $U$ admits a faithful, orientation-preserving action on a closed
4-manifold $M$ as in Theorem 4 then $q = 2$ or $q = 4$.}

\medskip

{\it Proof.}  i)  See [Z1, Proof of Prop.1].

ii) We denote by $g$ a generator of the normal subgroup $\Z_p$ of $U$.  If we are in
case ii) of Lemma 1 then $\Z_p$ fixes pointwise a 2-sphere $S^2$ which is invariant
under $\Z_q$ and $U$. If $q>2$ then $\Z_q$ and hence $U$ have a global fixed point
on  $S^2$; now a $U$-invariant regular neighbourhood in $M$ of this fixed point is a
3-sphere with a faithful action of $U$, contradicting part i) of the Lemma.

\medskip

If the fixed point set of $g$ consists of two 2-spheres then a subgroup of index two
of $U$ fixes pointwise both 2-spheres. As before, this is possible only for $q=2$ or
$q=4$.

\medskip

Finally, suppose that the fixed point set of $\Z_p$ consists of four isolated
points. Then a subgroup of index at most four of $U$ has a fixed point and acts
faithfully on a 3-sphere. Again by part i) of the Lemma, this is possible only for
$q$ = 2, 4 or 8. If $q=8$ then a dihedral subgroup $\Z_p \ltimes \Z_2$ of $U$ has a
fixed point in $M$. The situation for actions of dihedral groups has been analyzed
in [Mc1], in particular it follows from  [Mc1, Prop.13] that in the case $b_2(M)=2$
a dihedral group has to act without fixed points on $M$. This contradiction excludes
$q=8$ and finishes the proof of Lemma 3.

\bigskip

{\bf Lemma 4.} {\sl  i) If $G = {\rm PSL(2},p)$, for a prime $p\ge 5$, then
$p=5$ and $G$ is isomorphic to  $\A_5 \cong {\rm PSL(2,5)}$.

 ii) If $G = {\rm PSL(2},q)$, for
a prime power $q=p^n$ with $n>1$, then $q=4$ so again $G$ is isomorphic to
$\A_5 \cong  {\rm PSL(2,4)}$.

iii)  Let $\tilde G$ be a finite central extension, with nontrivial center, of
a nonabelian simple group $G$ (e.g., the central extension ${\rm SL(2},q)$ of
${\rm PSL(2},q))$. If $\tilde G$ acts faithfully on  $M$ then $G$ is isomorphic
to the dodecahedral group $\A_5$.}

\medskip

{\it Proof.}  i)  Since ${\rm PSL(2},p)$ has a metacyclic subgroup $U = \Z_p
\ltimes \Z_{(p-1)/2}$, with an effective action of $\Z_{(p-1)/2}$ on the normal
subgroup $\Z_p$, part i) follows from Lemma 3.
\medskip

ii) The  group $G = {\rm PSL}(2,p^n)$ has a subgroup $U = (\Z_p)^n \ltimes
\Z_{(q-1)/2}$ if $p$ is odd, resp. $U = (\Z_p)^n \ltimes \Z_{q-1}$ if $p=2$,
with an effective action of $\Z_{(q-1)/2}$ resp. $\Z_{q-1}$ on $(\Z_p)^n$.
Since ${\rm PSL(2},p)$ is a subgroup of ${\rm PSL}(2,p^n)$, part i) of the
Lemma implies that $p$ = 2, 3 or 5.
\medskip

Suppose  that $p$ = 3 or 5.  We consider the subgroup $(\Z_p)^n$ of $U$ and a
nontrivial element $g$ in $(\Z_p)^n$. If $g$ fixes pointwise a 2-sphere (cases ii)
and iii) of  Lemma 1) then this 2-sphere is invariant under $(\Z_p)^n$, hence there
is a faithful action of  $(\Z_p)^{n-1}$  on $S^2$ which is possible only if $n \le
2$. If $g$ fixes four isolated points then $(\Z_p)^n$ has a global fixed point, and
by Theorem 1 ii) again we have $n \le 2$.

\medskip

Suppose that $G =  {\rm PSL}(2,25)$, with a subgroup $U = (\Z_5)^2 \ltimes
\Z_{12}$. If the element $g$ in $(\Z_5)^2$ fixes pointwise one or two 2-spheres
then $(\Z_5)^2$ has two or four fixed points and some nonabelian subgroup of
$U$ has a global fixed point contradicting Theorem 1. If $g$ fixes four
isolated points then also $(\Z_5)^2$ fixes these points, and again a nonabelian
subgroup of $U$ has a global fixed point. So ${\rm PSL}(2,25)$ does not occur.

\medskip

Next we consider ${\rm PSL}(2,9)$, with a subgroup $(\Z_3)^2 \ltimes \Z_{4}$;
let $g$ be a nontrivial element in $(\Z_3)^2$. If $g$ fixes four isolated
points then all nontrivial elements in $(\Z_3)^2$ have exactly four fixed
points (since all subgroups $\Z_3$ are conjugate). The whole group $(\Z_3)^2$
fixes at least one of the four fixed  points of $g$ and hence admits a free
action on $S^3$; but $(\Z_3)^2$ does not admit a free action on $S^3$ (see [B])
so this case does not occur. A similar argument applies if $g$ fixes two
isolated points and a 2-sphere. Finally, if $g$ fixes pointwise two 2-spheres
then each of these 2-sphere is invariant under $(\Z_3)^2$ and $(\Z_3)^2$ has
two global fixed points on it. By [Mc1, Prop.14], the singular set of
$(\Z_3)^2$ consists of exactly four 2-spheres intersecting pairwise at their
poles; this contradicts the fact that all subgroups $\Z_3$ of $(\Z_3)^2$ are
conjugate and hence fix pointwise two 2-spheres. So ${\rm PSL}(2,9)$ does not
occur.

\medskip

This leaves us with the groups  ${\rm PSL}(2,2^n)$, for $n \ge 3$, with a
subgroup $U = (\Z_2)^n \ltimes \Z_{q-1}$ such that all involutions in
$(\Z_2)^n$ are conjugate. Suppose that $n \ge 4$. Let $g$ be an involution in
$(\Z_2)^n$.  If $g$ has four isolated fixed points or fixes two points and a
2-sphere, then a subgroup $(\Z_2)^2$ of $(\Z_2)^n$ has a global fixed point and
acts freely on $S^3$ which is a contradiction. Suppose that the fixed point set
of $g$ consists of two 2-spheres; then another involution in $(\Z_2)^n$ leaves
each of these 2-spheres invariant and acts  orientation-preservingly on it (by
[Mc1, Theorem 2]), any involution has to act orientation-preservingly  on such
a 2-sphere since the action is homologically trivial). Now again  a subgroup
$(\Z_2)^2$ has a global fixed point, and a contradiction to [Mc1, Prop.14] is
obtained as in the previous case of the subgroup $(\Z_3)^2$ of ${\rm
PSL}(2,9)$.

\medskip

Finally, we exclude the group ${\rm PSL}(2,8)$ which has a subgroup $(\Z_2)^3
\ltimes \Z_7$ such that all involutions are conjugate; let $g$ be an involution
in $(\Z_2)^3$. If $g$ fixes two 2-spheres then a subgroup $(\Z_2)^2$ has four
fixed points. By [Mc1, Prop.14] the singular set of $(\Z_2)^2$ is a union of
four 2-spheres which is a contradiction since each involution in $(\Z_2)^2$
fixes pointwise two 2-spheres.  If $g$ fixes a 2-sphere and two isolated points
then again a subgroup $(\Z_2)^2$ has four fixed points which is a contradiction
to [Mc1, Prop.14] since each involution fixes exactly one 2-sphere.
\medskip

Suppose that $g$ has four isolated fixed points. None of these points is fixed
by a subgroup $(\Z_2)^2$ since otherwise $(\Z_2)^2$ would act freely on a
3-sphere $S^3$. So each orbit under $(\Z_2)^3$ of a fixed point of an
involution has exactly four elements, and there are exactly seven such orbits.
This is exactly the situation excluded for a group $(\Z_2)^2$ in the proof of
[Mc2, Lemma 4.5], so ${\rm PSL}(2,8)$ does not act.

\medskip

iii) Let $g$ be a nontrivial central element of $\tilde G$. If $g$ fixes pointwise
one or two 2-spheres then such a 2-sphere is invariant under the factor group
$\tilde G/<g>$ or under a subgroup of index two, and hence $G \cong \A_5$. If $g$
has four isolated fixed points then  a subgroup of index at most four of $\tilde G$
fixes each of these four points; by Theorem 1, such a group has to be abelian so
this case does not occur.
\medskip

This finishes the proof of Lemma 4.

\bigskip
\bigskip

Continuing with the proof of Theorem 4 we consider next the groups  $G = {\rm
PSL}(m,q)$ in the Gorenstein-Harada list, where $q=p^n$ is odd or $m=3$ and
$q=4$. We note that ${\rm PSL}(m,p)$ is a subgroup of ${\rm PSL}(m,q)$; also,
for $r<m$, the linear group ${\rm SL}(r,q)$ is a subgroup of the linear
fractional group ${\rm PSL}(m,q)$ (see [Su2,chapter 6.5]). Applying Lemma 3 and
Lemma 4, it suffices then to exclude the groups ${\rm PSL}(3,4)$, ${\rm
PSL}(3,3)$ and ${\rm PSL}(3,5)$; but ${\rm PSL}(3,4)$ has a subgroup ${\rm
PSL}(3,2) \cong {\rm PSL}(2,7)$, the group ${\rm PSL}(3,3)$ has a metacyclic
subgroup $\Z_{13} \ltimes \Z_3$ and ${\rm PSL}(3,5)$ a metacyclic subgroup
$\Z_{31} \ltimes \Z_3$ which are all excluded by Lemma 3 or Lemma 4 (see [C]
for information about the subgroup structure  of the finite simple groups).
Thus among the linear fractional groups ${\rm PSL}(m,q)$ there remains only the
group ${\rm PSL}(2,5) \cong \A_5$.

\medskip

The proof for the unitary groups ${\rm PSU}(m,q)$ and the symplectic groups
${\rm PSp}(4,q)$ is similar, noting that ${\rm PSU}(2,q) \cong {\rm PSL}(2,q)
\cong {\rm PSp}(2,q)$. The unitary groups ${\rm PSU}(3,3)$ and ${\rm PSU}(3,5)$
are excluded since both have a subgroup ${\rm PSL}(2,7)$, and ${\rm PSU}(3,4)$
because it has a subgroup $\Z_{13} \ltimes \Z_3$. Noting that ${\rm SU}(r,q)$
is a subgroup of ${\rm PSU}(m,q)$, for $r<m$, by Lemma 4 iii) this excludes all
unitary groups ${\rm PSU}(m,q)$ except ${\rm PSU(2,5)} \cong  {\rm PSL}(2,5)$.

\medskip

Concerning the symplectic groups ${\rm PSp}(4,q)$, $q$ odd, we note that ${\rm
PSp}(4,3) \cong {\rm PSU}(4,2)$ has a subgroup $(\Z_3)^3$ and ${\rm PSp}(4,5)$
a subgroup $(\Z_5)^3$. Choosing an element $g$ in this subgroup $(\Z_3)^3$
resp. $(\Z_5)^3$ and applying Lemma 1, it is  easy to see that either such a
subgroup must have a global fixed point contradicting Theorem 1 ii), or there
is a subgroup $(\Z_3)^2$ resp. $(\Z_5)^2$ acting faithfully on a 2-sphere which
gives again a contradiction; hence these groups do not act on $M$. Since ${\rm
Sp}(2,q)$ is a subgroup of ${\rm PSp}(4,q)$ this excludes also the symplectic
groups ${\rm PSp}(4,q)$.

\medskip

Considering the remaining groups in the Gorenstein-Harada list, up to central
extension there are inclusions $^3{\rm D}_4(q) \supset {\rm G}_2(q) \supset
{\rm PSL}(3,q)$ which excludes these groups (see [St,Table 0A8], [GL,Table
4-1]). The Ree groups $^2{\rm G}_2(3^{2m+1})$ have one conjugacy class of
involutions, the centralizer of an involution is $\Z_2^2 \times {\rm
PSL}(2,3^{2m+1})$ ([G2,p.164]) so for $m \ge 1$ they do not act (the group
$^2{\rm G}_2(3)$ is not simple).

\medskip

The Sylow 2-subgroup $S_2$ of the Suzuki group Sz(8) has order 64 and a normal
subgroup $(\Z_2)^3$, and all involutions are conjugate. Let $g$ be a central
involution in $S_2$. If $g$ has four isolated fixed points, or fixes two isolated
points and a 2-sphere, then  a subgroup of order at least 16 fixes one of these
isolated fixed points and hence is abelian by Theorem 1. Since Sz(8) has no
elements of order eight, a subgroup $(\Z_2)^2$ fixes a point and hence admits a free
action on $S^3$ which is a contradiction. If the fixed point set of $g$ consists of
two 2-spheres then a subgroup of order at least 16 of $S_2$ acts faithfully and
orientation-preservingly on a 2-sphere; since there are no elements of order eight
this gives again a contradiction, so Sz(8) does not occur.

\medskip

Finally, $\A_7$ has a subgroup ${\rm PSL}(2,7)$,  and also the Mathieu groups
${\rm M}_i$, the Janko groups ${\rm J}_i$, the McLaughlin group McL and the
Lyons group Ly have metacyclic or linear fractional subgroups excluded by Lemma
3 or Lemma 4 (see [C]).

\bigskip

Hence we have excluded all finite simple groups from the Gorenstein-Harada list
except the alternating group  $\A_5$,  and for the proof of Theorem 4 it remains to
show that the 4-manifold $M$ has intersection form of $S^2 \times S^2$.

\bigskip

Now $G = \A_5$ has a subgroup $\A_4  \cong  (\Z_2 \times \Z_2) \ltimes \Z_3$, and we
consider the normal subgroup $(\Z_2)^2$ of $\A_4$. By [Mc1, Prop.14] either $M$ has
the right intersection form or $(\Z_2)^2$ has a global fixed point, so we can assume
the latter.  Let $g$ be an involution in $(\Z_2)^2$. According to Lemma  1 we
consider three cases.

\medskip

If $g$ has four isolated fixed points then $(\Z_2)^2$ has two or four global fixed
points. But then also $\A_4$ has a global fixed point which contradicts Theorem 1.

\medskip

If $g$ fixes pointwise two 2-spheres then each of these 2-spheres is invariant under
the action of $(\Z_2)^2$; moreover by [Mc1, Theorem 2], since $(\Z_2)^2$ acts
homologically trivial it acts orientation-preservingly on each of these 2-spheres.
Then $(\Z_2)^2$ has again exactly four global fixed points and we get a
contradiction as in the first case.

\medskip

Finally, if the fixed point set of $g$ consists of a 2-sphere and two isolated
points then $(\Z_2)^2$ has two or four global fixed points, so $\A_4$ has a global
fixed point contradicting Theorem 1.

\bigskip

{\bf 3. Proof of Theorem 5}

\medskip

Recall that a finite $Q$ group is {\it quasisimple} if it is perfect (the
abelianized group is trivial) and the factor group of $Q$ by its center is a
nonabelian simple group.  A finite group $E$ is  {\it semisimple} if it is perfect
and the factor group of $E$ by its center is a direct product of nonabelian simple
groups (see [Su2, chapter 6.6]). A semisimple group is a central  product of
quasisimple groups that are uniquely determined. Any finite group $G$ contains a
unique maximal semisimple normal group $E(G)$ (the subgroup $E(G)$ may be trivial);
the subgroup $E(G)$ is characteristic in $G$ and the quasisimple factors of $E(G)$
are called the {\it components} of $G$.     To prove the Theorem 5  we consider
first the case of  trivial maximal normal semisimple subgroup and prove that in
this case the groups are solvable.

\bigskip

{\bf  Lemma 5.}  {\sl Let $G$ be a finite group with  trivial maximal normal
semisimple subgroup $E(G)$. If $G$  admits  a  homologically trivial action on
a closed 4-manifold M with $b_2(M)=2$ and trivial first homology, then $G$ is
solvable.}

\medskip

{\it Proof.}   We consider first the  case of $G$ containing  a normal nontrivial
cyclic subgroup $H$ and we prove that in this case if $E(G)$ is trivial then $G$ is
solvable. We can suppose that
$H$ has prime order $p$ so each nontrivial element of $H$ has the same fixed point
set; since $G$ normalizes $H$ then $G$  fixes setwise the fixed point set of $H$.

\medskip

If the fixed point set of $H$ consists of  two isolated points and a 2-sphere there
exists a subgroup  $G_0$  of index at most two in $G$ such that $G_0$ fixes both
points; the subgroup $G_0$ is abelian by Theorem 1 and consequently $G$ is solvable.

\medskip

If the fixed point set of $H$ consists of  four  isolated points there exists a
normal subgroup $G_0$ of $G$   that fixes each  point; $G_0$ is abelian by
Theorem 1. The quotient group $G/G_0$ is isomorphic to a subgroup of
$\Bbb{S}_4$, the symmetry group over four elements that is a solvable group and
we can conclude that $G$ is solvable. Finally suppose that the fixed point set
of   $H$ consists of  two 2-spheres;  there exists a subgroup  $G_0$  of index
at most two in $G$ such that $G_0$  leaves invariant both 2-spheres. We
consider in $G_0$  the normal subgroup $K$ of elements fixing pointwise
$S^2_+$, one of the two 2-spheres; the subgroup $K$ contains $H$ and since $K$
acts locally by rotations around $S^2_+$ then $K$ is cyclic. The factor group
$G_0/K$ acts faithfully on $S^2_+$. If $G_0/K$ is solvable, we get the thesis;
otherwise we can suppose that  $G_0/K$  is isomorphic  to  $\Bbb{A}_5$ because
it is  the only nonsolvable finite group acting orientation preservingly on the
2-sphere (the action is orientation preserving by a result of Edmonds  [E], see
also  [Mc1,  Theorem 2]).       The action of $\Bbb{A}_5$ by conjugation  on
$K$ is trivial because $K$ is cyclic and its automorphism group is abelian;
then $G_0$ is a central extension of $\Bbb{A}_5$. The derived group $G'_0$ is a
quasisimple normal subgroup of $G_0$  (see [Su1, Theorem 9.18, pag.257]); this
fact implies that  $E(G_0)$,  and consequently $E(G)$, are not trivial  in
contradiction with our hypothesis.   The proof of this particular case is now
complete  and in the following we can use this fact.

\medskip

{\it Fact: if a subgroup $N$ of $G$ contains  a nontrivial cyclic normal subgroup
then either  $N$ is solvable or  $E(N)$ is not trivial.}

\medskip

We consider now the general case.  We denote by $F$  the Fitting subgroup of $G$
(the maximal nilpotent normal subgroup of $G$). Since  $E(G)$ is trivial,  the
Fitting subgroup
$F$ coincides with   the generalized Fitting subgroup that  is the product of the
Fitting subgroup with the maximal semisimple normal subgroup. The generalized
Fitting subgroup $F$ contains its centralizer in $G$  and
$F$ is not trivial ([Su2, Theorem 6.11, pag.452]).

\medskip

Since $F$ is nilpotent it is the direct product of its Sylow $p$-subgroups. In
particular any Sylow subgroup of $F$ is normal in $G$; since $F$ is not trivial
we have  a nontrivial $p$-subgroup $P$ which is normal in $G$. We consider the
maximal elementary abelian $p$-subgroup $Z$ contained in the center of $P$;
this subgroup is not trivial and it is normal in $G$.

\medskip

Suppose first that we can chose $p$  odd (the order of $F$ is not a power of
two). If $Z$ contains an element with fixed point set consisting  of four
points or of two points and one 2-sphere, the group $Z$ has global fixed point
set and has  rank at most two by Theorem 1.  If $Z$ contains an element with
fixed point set consisting of  two 2-spheres each element of $Z$ leaves
invariant both 2-spheres; a quotient group of $Z$ by a cyclic group acts
faithfully on the 2-spheres. This quotient group has to be  cyclic and it acts
on the 2-spheres by rotations then the group $Z$ has also in this case global
fixed point set and it has rank at most two. If $Z$ is cyclic, by the first
part of the proof,   $G$ is solvable. If $Z$ has rank two, since it has global
fixed point set it is described by [Mc1, Prop.14]. The fixed point set of $Z$
consists of four points. The whole group $G$ leaves invariant the fixed point
set of $Z$    and there exists a normal subgroup $G_0$ that fixes each point.
The quotient group $G/G_0$ is isomorphic to a subgroup of $\S_4$ that is
solvable. By Theorem 1  $G_0$ is abelian  and consequently $G$ is  solvable.

\medskip

Suppose now that the order of $F$ is a power of two; in this case $F=P$ is a
2-group and $Z$ is an elementary abelian 2-group of  rank at most four (by
Lemma 2).  If $Z$ has rank one by the first part of the proof $G$ is solvable.
If $Z$ has rank two we consider $C_G(Z)$ the centralizer of $Z$ in $G$ that is
normal because $Z$ is normal;   $C_G(Z)$ contains a nontrivial normal cyclic
subgroup and it is solvable. The factor $G/C_G(Z)$ is isomorphic to  a subgroup
of GL$(2,\Bbb{Z}_2)$, the automorphism group of an elementary abelian 2-group
of rank two; since GL$(2,\Bbb{Z}_2)$ is a solvable group we can conclude that
$G$ is solvable.  Suppose that  $Z$ has rank three.  In this case the factor
group $G/C_G(Z)$ is  isomorphic to  a subgroup of GL$(3,\Bbb{Z}_2)$, the
automorphism group of an elementary abelian 2-group of rank three;
GL$(3,\Bbb{Z}_2)$ has order $2^3\cdot 3\cdot 7$ and any element of order seven
permutes cyclically all the involutions of $(\Bbb{Z}_2)^3$. The group
$G/C_G(Z)$ can not contain any element of order $7$ otherwise all involutions
in $Z$  are conjugated and this can be excluded by the same argument used to
exclude case ${\rm PSL}(2,8)$ in the proof of Lemma 4; so the group $G/C_G(Z)$
has order $2^\alpha 3^\beta$ and it is solvable. This fact implies that  $G$ is
solvable.

\medskip

It remains the case $Z$ of rank four; the factor group
$G/C_G(Z)$ is  isomorphic to  a subgroup of   GL$(3,\Bbb{Z}_2)$, the automorphism
group of an elementary abelian 2-group of rank four.

\medskip

We analyze the fixed point set of the elements in $Z$. The group   $Z$ cannot
contain any element with fixed point set consisting  of  two points and one
2-sphere. In this case a subgroup of index at most two of $Z$ has global fixed
point set and this impossible by Theorem 1.

\medskip

Suppose first that   $h$ is  an involution such that its fixed point set Fix($h$)
consists of two 2-spheres.   If an element of $Z$ leaves invariant both components
of Fix($h$), it acts on both 2-spheres orientation preservingly (by [Mc1,  Theorem
2]). Since there exists only one involution acting trivially on Fix($h$) and the
maximal elementary 2-group acting faithfully on a 2-sphere has rank two, the group
$Z$ contains with index two a subgroup $Z_0$ that leaves invariant both 2-spheres in
Fix($h$).  Any involution of $Z_0$ different from $h$ acts  nontrivially and
orientation preservingly (again by  [Mc1, Theorem 2]) on the 2-spheres, thus it
fixes pointwise on each 2-sphere two points, this implies that the subgroup of rank
two generated by $h$ and by the other involution  has global fixed point set; this
2-rank subgroup is described by [Mc1, Prop.14] and it contains two involutions
different from $h$, one with 0-dimensional fixed point set and one with
2-dimensional fixed point set.  We obtain that if the  fixed point set of $h$
consists of two 2-spheres there exist exactly three involutions with  0-dimensional
fixed point set and with fixed point set contained in Fix($h$).   We consider now
$h'$ an involution such that its fixed point set Fix($h'$) consists of four isolated
points. We consider $Z_1$ the subgroup of $Z$ that fixes pointwise Fix($h'$);  since
a maximal elementary 2-group in the symmetric group  $\Bbb{S}_4$  has rank two, the
subgroup  $Z_1$ has index at most four but by Theorem 1 the subgroup $Z_1$ has rank
at most two; we can conclude that the rank of $Z_1$ is exactly two. In this case
$Z_1$ is completely described by  [Mc1, Prop. 14] and there exist exactly two
involutions in $Z$ with 2-dimensional fixed point set that contain Fix($h'$).   If
$n$ is the number of involutions in $Z$ with 0-dimensional fixed point set and $m$
is the number of involutions with 2-dimensional fixed point set we have that
$n+m=15$ and for the previous computations $3m=2n$; we obtain that $n=9$ and
$m=6$.   We recall that the factor group
$G/C_G(Z)$ is  isomorphic to  a subgroup of GL$(3,\Bbb{Z}_2)$; the group
GL$(3,\Bbb{Z}_2)$  has order
$2^6\cdot 3^2 \cdot 5\cdot 7$;  an automorphism of order five does not centralize
any involution of $(\Bbb{Z}_2)^3$ (we have three orbits with five elements) and an
automorphism of order seven centralizes exactly one involution (two orbits with
seven elements and one orbit with only one element).  Since we have nine  elements
with 0-dimensional fixed point set and six elements with 2-dimensional fixed point
set, the group  $G/C_G(Z)$ can not contain elements of order five and seven and it
has order $2^\alpha 3^\beta$; we have that   $G/C_G(Z)$ is solvable and consequently
$G$ is solvable. This fact concludes the proof.

\bigskip

The following lemma considers  the case of semisimple groups.

\bigskip

{\bf  Lemma 6.}  {\sl Let $G$ be a finite semisimple group that   admits  a
homologically trivial action on a closed 4-manifold M with $b_2(M)=2$ and
trivial first homology; then $G$ is  isomorphic to one of the groups $\;
\Bbb{A}_5,\; \Bbb{A}_5^*$ or $\; \Bbb{A}_5\times \Bbb{A}_5.$ }

\medskip

{\it Proof.} By Lemma 4(iii),   if $G$ is quasisimple, then $G$  is isomorphic
either to $\Bbb{A}_5$  or to
$\Bbb{A}_5^*\cong$SL(2,5) that is the unique perfect central extension of
$\Bbb{A}_5$.

\medskip

We consider now the case of   $G$ with two quasisimple components; since in our list
of quasisimple groups
$\Bbb{A}_5^*$  is the unique group with nontrivial center, then either $G\cong
\Bbb{A}_5^*\times_{\Bbb{Z}_2}\Bbb{A}_5^*$ or $G$ is the direct product of two
quasisimple subgroups.

\medskip

We prove  that  no involution can be contained in the center of  $G$   and the only
possibility with two components remains $\Bbb{A}_5\times \Bbb{A}_5$.   Suppose that
$h$ is  an involution contained in the center of $G$. If the fixed point set
Fix($h$) of $h$ consists of two points and one 2-sphere then   $G$ has a subgroup of
index at most two  that fixes both points and by Theorem 1 this group should be
abelian of rank two, that is impossible.

\medskip

 If the fixed point set of $h$ consists of four isolated points there exists a
normal subgroup $G_0$ of $G$ that  fixes each point of Fix($h$); by Theorem 1
$G_0$ is abelian of rank two. The quotient $G/G_0$ acts faithfully on the four
points of Fix($h$) and it is isomorphic to a subgroup of $\Bbb{S}_4$. The group
$G$ contains a normal subgroup isomorphic either to $\Bbb{A}_5^*$ or to
$\Bbb{A}_5$; no quotient of  these groups  by an abelian normal subgroup is
isomorphic to a subgroup of $\Bbb{S}_4$. Finally we consider when Fix($h$)
consists of two 2-spheres, in this case $G$  leaves invariant both 2-spheres
because it does not contain any subgroup of index two ($G$ is perfect). The
subgroup acting trivially on the 2-spheres is cyclic and normal in $G$; the
quotient of $G$ by it acts faithfully on the 2-spheres and it is again the
product of two quasisimple groups. This fact can not occur.  By the previous
part if $G$ has three or more components the quasisimple factors are all
isomorphic to $\Bbb{A}_5$ but these groups cannot occur because the sectional
2-rank of $G$ is smaller then four.

\bigskip

{\it Proof of Theorem 5.}

\medskip

If the maximal semisimple normal subgroup $E(G)$ is trivial then $G$ is solvable by
Lemma 5. We can suppose that $E(G)$ is not trivial and it is isomorphic by Lemma 6
to  $\Bbb{A}_5,\, \Bbb{A}_5^*$ or $\Bbb{A}_5\times \Bbb{A}_5$.  We denote by $C$ the
centralizer of $E(G)$ in $G$; since $E(G)$ is normal its centralizer is normal.
Whatever is $G$  the fixed point sets of the elements in the centralizer consist of
two 2-spheres. Let $h$ be a nontrivial  element in $C$. If the fixed point set of
$h$ consists of two points and one 2-sphere then   $G$ has a subgroup of index at
most two  that fixes both points and by Theorem 1 this group should be abelian of
rank two, that is impossible.
 If the fixed point set of $h$ consists of four isolated points there exists a
normal subgroup $G_0$ of $G$ that  fixes each point of Fix($h$); by Theorem 1 the
group $G_0$ is abelian of rank two. The quotient $G/G_0$ acts faithfully on the four
points of Fix($h$) and it is isomorphic to a subgroup of $\Bbb{S}_4$. The group $G$
contains a normal subgroup isomorphic either to $\Bbb{A}_5^*$ or to $\Bbb{A}_5$; no
quotient of  these groups  by an abelian normal subgroup is isomorphic to a subgroup
of $\Bbb{S}_4$.

\bigskip

{\it Case 1: $E(G)$ is isomorphic to $\Bbb{A}_5\times \Bbb{A}_5$}

\medskip

If   $E(G)$ is isomorphic to $\Bbb{A}_5\times \Bbb{A}_5$ the subgroup $C$ is
trivial. In fact if $C$ contains a nontrivial element $f$ all the group
$\Bbb{A}_5\times \Bbb{A}_5$ fixes setwise the fixed point set of $f$. The group
$E(G)$ does  not contain any subgroup of index two (it is perfect), so
$\Bbb{A}_5\times \Bbb{A}_5$ fixes setwise each 2-sphere; the group  $E(G)$ does not
contain any nontrivial normal cyclic subgroup so it acts faithfully on  both
2-spheres and this is impossible. This fact implies that $G$ is isomorphic to a
subgroup of  Aut$(E(G))$, the automorphism group of $E(G)$,that contains with index
two the subgroup  Aut($\Bbb{A}_5$)$\times$  Aut($\Bbb{A}_5$)$\cong \Bbb{S}_5\times
\Bbb{S}_5$ and any element not in $\Bbb{S}_5\times \Bbb{S}_5$ exchange the two
quasisimple components of $E(G)$ (see  [GLS, Theorem 3.23, p.13]). We recall that
the quasisimple components of a group are permuted by any automorphism of the group
(see  [GLS, Theorem 3.5, p.7]).   We consider $S$  the Sylow 2-subgroup of $E(G)$
that is an elementary abelian 2-group of rank two; in the proof of Lemma 5 we have
proved that $S$ contains six elements with fixed point set consisting  of two
2-spheres  and nine elements with  fixed point set consisting of four points.  Let
$h$ be an  element with 2-dimensional fixed point set; $h$  is contained in a
subgroup $A$ of rank two with global fixed point set. By  [Mc1, Prop. 14] each of
the four 2-spheres in the singular set of $A$ represents a primitive class in
$H_2(M)$ and together these classes generate $H_2(M)$. Since the 2-spheres of the
fixed point of $h$ are exchanged by some elements of $S$ and the action is
homologically trivial we can conclude that the two 2-spheres of the fixed point of
$h$ represent the same class. Each element with 2-dimensional fixed point set gives
one primitive class in homology and the 2-spheres in the singular set of $S$
generate  $H_2(M)$ (that has rank two).   The conjugacy classes of involutions of
$S$ in $E(G)\cong \Bbb{A}_5\times \Bbb{A}_5$ are three. The six involutions
contained in the  quasisimple components form two conjugacy classes, one class for
each component; the third class contains the remaining nine involutions.  Since the
action is homologically trivial, we have that the six involutions in the quasisimple
components have 2-dimensional fixed point set and the two  conjugacy classes in the
quasisimple components represent  different elements in  $H_2(M)$. This fact implies
immediately  that $G$ cannot contain any element that by conjugation exchanges  the
two quasisimple components and $G$ is isomorphic to  a subgroup of
Aut($\Bbb{A}_5$)$\times$  Aut($\Bbb{A}_5$)$\cong \Bbb{S}_5\times \Bbb{S}_5$.
Moreover  we consider $t$ an involution contained in a quasisimple component.  The
fixed point set of $t$, Fix($t$), consists of two 2-spheres.  The other involutions
contained in $S$ that are in the same  quasisimple component of $t$ act on Fix($t$)
exchanging the two 2-spheres otherwise the Sylow 2-subgroup of the quasisimple
component has global fixed point set and it  contains an involution with
0-dimensional fixed point set that is impossible.  Suppose that  a subgroup
isomorphic to $ \Bbb{S}_5\times \Bbb{A}_5$ is contained in $G$. We consider $t$ an
involution contained in the second factor (that  is  a quasisimple component of the
group), its centralizer in $G$ contains a subgroup isomorphic to  $ \Bbb{S}_5\times
\Bbb{Z}_2\times \Bbb{Z}_2$. Since the involutions in the same quasisimple component
exchange the two 2-spheres of Fix($t$) there exists a group isomorphic to  $
\Bbb{S}_5\times \Bbb{Z}_2$ that leaves invariant the two 2-spheres of Fix($t$).  The
normal subgroup of  $ \Bbb{S}_5\times \Bbb{Z}_2$ that acts trivially on the
2-spheres has to be cyclic;  its intersection with the first factor has to be
trivial.  We have that $\Bbb{S}_5$ acts faithfully on the two 2-spheres and this is
impossible.  We obtain that either $G\cong \Bbb{A}_5\times \Bbb{A}_5$ or $G$
contains with index two $E(G)$ and $G$ is isomorphic to the  extension of $
\Bbb{A}_5\times \Bbb{A}_5$ by the outer automorphism that  acts nontrivially on
both quasisimple components.
\bigskip

{\it Case 2: $E(G)$ is isomorphic to $\Bbb{A}_5$}

\medskip

We consider $F(G)$ the Fitting subgroup of $G$; the Fitting subgroup centralizes
$E(G)$ [Su2,  pag.452], then all the elements in the Fitting have two 2-spheres as
fixed point set.   Suppose  first that $F(G)$ contains a nontrivial cyclic
subgroup  normal in $G$; in this case all the group fixes setwise the two 2-spheres
that are the fixed point set of this group. We consider $K$ the normal subgroup
consisting of the elements acting trivially on the two 2-spheres; the subgroup $K$
is cyclic. The subgroup  $E(G)\cong \Bbb{A}_5$ intersects trivially $K$ and since
$E(G)$   does not contain subgroup of index two it leaves invariant  both
2-spheres.  We recall that by [Mc1, Theorem 2] the action of $G$ on the two
2-spheres  has to be orientation preserving and  $\Bbb{A}_5$  is maximal between the
finite groups acting orientation preserving on a  2-sphere. Moreover since the
automorphism group of $K$ is abelian, the group $E(G)$ acts trivially by conjugation
on $K$. We can conclude in this case that the subgroup of $G$ of index at most two
that leaves invariant both the 2-spheres in the fixed point of $K$  is exactly
$E(G)\times K$.   We supposse now that the order of the Fitting subgroup is not a
power of two; the Fitting subgroup contains a characteristic $p$-subgroup with $p$
odd. The center of this subgroup is normal in $G$ and it has to be cyclic, otherwise
the centralizer of  any element in the center contains a subgroup isomorphic to
$\Bbb{A}_5\times \Bbb{Z}_{p}\times\Bbb{Z}_{p}$ that is impossible because
$\Bbb{A}_5$ is maximal between the finite  groups acting orientation preservingly on
a 2-sphere. We are in the case considered previously and we get the thesis.  We can
suppose finally that the order of $F(G)$ is a power of two. We consider $Z$ the
center of $F(G)$ that is a normal subgroup of  $G$. Since the sectional 2-rank of
the group $G$ is at most four by Lemma 2 and $E(G)$ contains an elementary
2-subgroup of rank two, we obtain that $Z$ has rank at most two. If it is cyclic we
have the thesis, so we can suppose that $Z$ has rank two.       We consider the
fixed point set of an involution in $Z$ that consists of two 2-spheres;  by the same
argument used previously we can obtain that the subgroup of index  at most  two  of
$E(G)\times F(G)$ that leaves invariant both the 2-spheres is isomorphic to
$\Bbb{A}_{5}\times\Bbb{Z}_{2^m}$. Since $Z$ has rank two we obtain that
$F(G)\cong\Bbb{Z}_{2^m}\times \Bbb{Z}_{2}$. If $m>1$ there exists a central
involution  in $G$ (in $F(G)$ the exist a unique involution that is not primitive)
and we get the thesis. It remains the case when $F(G)\cong \Bbb{Z}_{2}\times
\Bbb{Z}_{2}$. In this case $F(G)$ is the center of the generalized Fitting subgroup
$F^*(G)$ [Su2, Theorem 6.10., pag.452], $G/F(G)$ is isomorphic to a subgroup of
Aut$(F^*(G))\cong \Bbb{S}_5\times$ GL$(2,\Bbb{Z}_2)\cong \Bbb{S}_5 \times
\Bbb{S}_3$.
 We can suppose in $G/F(G)$ the existence of  an  element of order three   permuting
the three involutions in $F(G)$, otherwise one involution in $F(G)$ is centralized
by the whole group and we can conclude that $G$ contains a subgroup of index two
isomorphic to $\Bbb{A}_5\times \Bbb{Z}_{m}$. We denote by $\bar f \in G/F(G)$ such
an element  that we can suppose of order three. We denote by $f\in G$ a preimage of
$\bar f$ with respect to the projection of $G$ onto $G/F(G)$; we can chose $f$ of
order three.  Since the automorphisms of $ \Bbb{A}_5$ of order three are inner, we
can compose $f$ with the appropriate element of  $ \Bbb{A}_5$ and we can suppose
that the  action of $f$ on $E(G)\cong \Bbb{A}_5$ is trivial.   Now we want to prove
that the whole group Aut$(F^*(G))$ cannot occur as  the quotient  $G/F(G)$;  in
particular we exclude the presence in $G/F(G)$ of an involution acting  on
$\Bbb{A}_5$ as a non-inner  automorphism,  acting trivially on $F(G)$ and acting
trivially by conjugation on Aut$(F(G))$. Let $\bar t$ be an involution of this
type;  we denote by $t\in G $  a preimage of  $\bar t$ with respect to the
projection of $G$ onto $G/F(G)$. The element $t$ is  an involution  otherwise its
square is an involution in $F(G)$ and this is in  contradiction with the existence
of $f$ (we should have that $ft^2f^{-1}$ is an involution  different from $t^2$  and
on the other hand  $ftf^{-1}=t$). The element $t$ is an involution that acts
trivially by conjugation on $F(G)$, thus we have a subgroup of $G$ isomorphic to
$\S_5\times F(G)\cong \S_5\times \Bbb{Z}_{2}\times \Bbb{Z}_{2}$ where the fixed
point set of each  involution in $F(G)$ consists of two 2-spheres. In the Case 1 we
proved that  this group cannot occur.

We have that $E(G)\times F(G)$ is contained in $G$ with index at most six. The
split extension of $E(G)\times F(G)$ by $f$ is isomorphic to $\Bbb{A}_5\times
\Bbb{A}_4$ and it is contained in $G$ with index at most two.

\bigskip

{\it Case 3: $E(G)$ is isomorphic to $\Bbb{A}_5^*$}

\medskip We denote by $h$ the involution in the center of  $E(G)$ and by Fix($h$)
the fixed point set of $h$ that consists  consists of two 2-spheres.   The center of
$E(G)$ is normal  in $G$ and $G$ fixes setwise Fix($h$). We consider $G_0$ the
subgroup of $G$ of the elements that fixes setwise each 2-sphere in Fix($h$). The
subgroup $G_0$ has index at most two in $G$ and it contains $E(G)$.  In $G_0$ we
consider $K$ the cyclic and normal subgroup of the elements that act trivially on
one of the two 2-spheres; the quotient $G_0/K$ acts faithfully on the 2-sphere and
it contains a subgroup isomorphic to $\Bbb{A}_5$ (the quotient of $E(G)$ by its
center).  We recall that by [Mc1, Theorem 2] the action of $G_0$ on the two
2-spheres  has to be orientation preserving and  $\Bbb{A}_5$  is maximal between the
finite groups acting orientation preserving on a  2-sphere. We obtain that
$G_0/K\cong \Bbb{A}_5$. The group  $G_0$ is generated by $E(G)$ and $K$. The action
by conjugation of  $E(G)$ on $K$ has to be trivial because the automorphism group of
$K$ is abelian. This implies that $G_0$ is the central product of $F(G)$ and
$K$.

\bigskip \bigskip

\centerline {\bf References}

\bigskip

\item {[B]} G.Bredon, {\it Introduction to Compact Transformation Groups.} Academic
Press, New York 1972

\item {[C]} J.H.Conway, R.T.Curtis, S.P.Norton, R.A.Parker, R.A.Wilson, {\it Atlas
of Finite Groups.} Oxford University Press 1985

\item {[E]} A.L.Edmonds, {\it Aspects of group actions on four-manifolds.} Top.
Appl. 31, 109-124 (1989)

\item {[F]} M.H.Freedman, {\it The topology of four-dimensional manifolds.} J. Diff.
Geom. 17,  357-453  (1982)

\item {[G1]} D.Gorenstein, {\it The Classification of Finite Simple Groups.} Plenum
Press, New York 1983

\item {[G2]} D.Gorenstein, {\it Finite Simple Groups: An Introduction to Their
Classification.}  Plenum Press, New York 1982

\item {[GL]} D.Gorenstein, R.Lyons, {\it The local structure of finite group of
characteristic 2 type.} Memoirs Amer. Math. Soc. 276 (vol. 42), 1-731 (1983)

\item{[GLS]}   D. Gorenstein, R. Lyons, R. Solomon, {\it The classification of the
finite simple groups, Number 2.} Math. Surveys Monogr. 40, 2 (Amer. Math. Soc.,
Providence, RI, 1996)

\item {[HL]} I.Hambleton, R.Lee, {\it Finite group actions on $P^2(\Bbb C)$.}  J.
Algebra 116, 227-242  (1988)

\item {[Mc1]} M.P.McCooey, {\it Symmetry groups of 4-manifolds.} Topology 41,
835-851 (2002)

\item {[Mc2]} M.P.McCooey, {\it Groups which act pseudofreely on $S^2 \times S^2$.}
Pacific J. Math. 230, 381-408 (2007)  (arXiv:math.GT/9906159)

\item {[MZ1]} M.Mecchia, B.Zimmermann, {\it On finite simple and nonsolvable groups
acting on homology 4-spheres.} Top. Appl. 153,  2933-2942  (2006)

\item {[MZ2]} M.Mecchia, B.Zimmermann, {\it On finite groups acting on
$\Z_2$-homology 3-spheres.} Math. Z. 248, 675-693 (2004)

\item {[MZ3]} M.Mecchia, B.Zimmermann, {\it On finite simple groups acting
on integer and mod 2 homology 3-spheres.}  J. Algebra 298, 460-467  (2006)

\item {[St]} E.Stensholt, {\it Certain embeddings among finite groups of Lie type.}
J. Algebra 53, 136-187  (1978)

\item {[Su1]} M.Suzuki, {\it Group Theory I.}  Springer-Verlag 1982

\item {[Su2]} M.Suzuki, {\it Group Theory II.}  Springer-Verlag 1982

\item {[TD]} T. tom Dieck, {\it Transformation Groups.} Walter de Gruyter, Berlin
1987

\item {[W1]} D.M.Wilczy\'nski, {\it Group actions on the complex projective space.}
Trans. Amer. Math. Soc. 303,  707-731  (1987)

\item {[W2]} D.M.Wilczy\' nski, {\it Symmetries of homology complex projective
planes.}  Math. Z. 203, 309-319  (1990)

\item {[Z1]} B.Zimmermann, {\it On finite simple groups acting on homology
3-spheres.}  Top. Appl. 125, 199-202 (2002)

\item {[Z2]} B.Zimmermann, {\it On the classification of finite groups acting on
homology 3-spheres.}  Pacific J. Math. 217, 387-395 (2004)

\bye